\documentclass[final,1p,times]{elsarticle}
\usepackage{amssymb,color,ulem,amsmath}

\usepackage[usenames,dvipsnames]{xcolor}

\def\qed{\hfill $\square$}

\journal{}
\usepackage{amsthm,verbatim}
\theoremstyle{definition}

\newtheorem{thm}{Theorem}[section]

\newproof{pf}{Proof}
\newcommand{\eps}{\varepsilon}
\newcommand{\R}{{\mathbb{R}}}

\begin{document}
\begin{frontmatter}

\author[mymainaddress1]{Oleg Makarenkov\corref{mycorrespondingauthor}}
\cortext[mycorrespondingauthor]{Corresponding author}
\ead{makarenkov@utdallas.edu}


\address[mymainaddress1]{Department of Mathematical Sciences, University of Texas at Dallas, 75080, TX, Richardson, USA}


\title{Bifurcation of limit cycles from a switched equilibrium in planar switched systems}
\begin{abstract} We consider a switched system of two subsystems that are activated as the trajectory enters the regions $\{(x,y):x>\bar x\}$ and $\{(x,y):x<-\bar x\}$ respectively, where $\bar x$ is a positive parameter. We prove that a regular asymptotically stable equilibrium of the associated Filippov equation of sliding motion (corresponding to $\bar x=0$) yields an orbitally stable limit cycle for all $\bar x>0$ sufficiently small. Such an equilibrium is called   switched equilibrium in control theory, in which case  considering $\bar x>0$ refers to the effect of hysteresis. The fact that hysteresis perturbation of a switched equilibrium yields a limit cycle is known and is actively used in control. We not only prove this fact rigorously for the first time ever, but also offer a formula for the period of the limit cycle which turns out to be very sharp as our practical example  demonstrates. Specifically, an application to a model of a dc-dc power converter concludes the paper.
\end{abstract}
\end{frontmatter}
\section{Introduction}\label{sec:int}

\begin{figure}[h]\center
\includegraphics[scale=0.9]{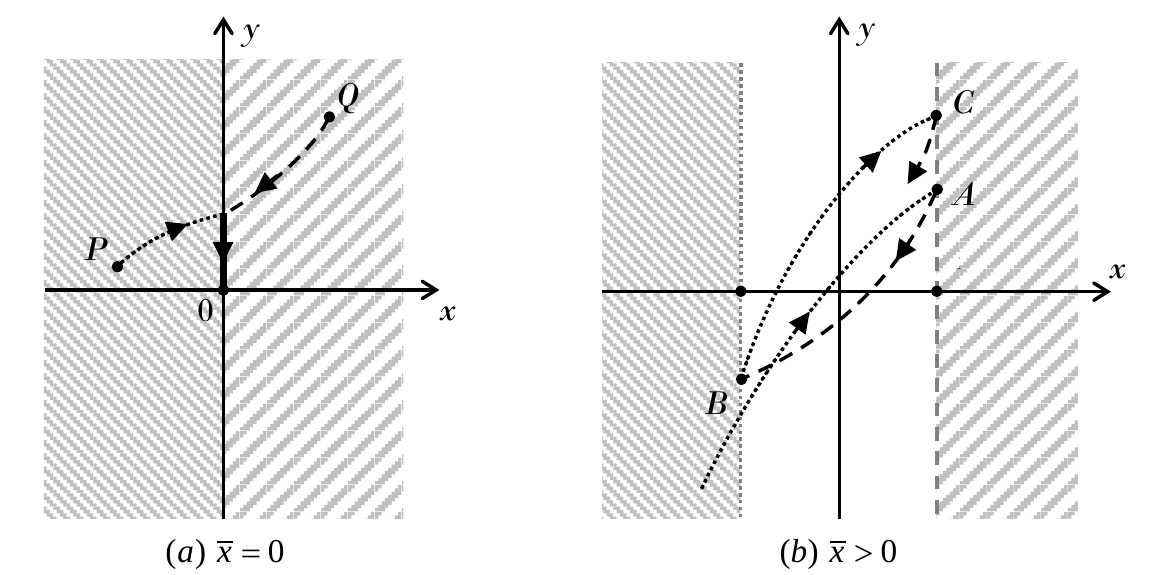}
\vskip-0.2cm
\caption{\footnotesize  Trajectories of switched system (\ref{fg})-(\ref{RL}) for (a) $\bar x=0$, when the trajectories with any initial condition (e.g. $P$ and $Q$ as shown) stick to the switching threshold and slip along it until reaching the origin, and for (b) $\bar x>0$, when the trajectory switches to $k=+1$  ($k=-1$) upon reaching $x=\bar x$ ($x=-\bar x$). The small and coarse textures stay for the regions of  $k=-1$ and $k=+1$ respectively.} \label{fig0}
\end{figure}

\noindent  The paper investigates the existence of attracting limit cycles in switched systems of the form
\begin{eqnarray}
&&  \begin{array}{l}
     \dot x(t)=f^k(x(t),y(t)),\\
    \dot y(t)=g^k(x(t),y(t)),
\end{array} \label{fg}\\
&&\begin{array}{l}
     k:=+1,\quad{\rm if }\ \  x(t)>\bar x,\\
     k:=-1,\quad{\rm if }\ \ x(t)<-\bar x,
  \end{array}\label{RL}
\end{eqnarray}
where $f^-,$ $f^+,$ $g^-,$ $g^+$ are smooth functions and $\bar x\in\mathbb{R}$ is a parameter. The formulas (\ref{fg})-(\ref{RL}) is a shorthand for a switching rule that operates as follows. To start drawing a trajectory of system (\ref{fg})-(\ref{RL}) one  needs the initial point $(x(0),y(0))$ and  the index of the system ($k=-1$ or $k=+1$) that governs the trajectory at $t=0.$ The trajectory then follows  system $i$ until it reaches one of the lines $\{-\bar x\}\times\R$ or $\{\bar x\}\times\R$, when $k$ switches to $k=-1$ or $k=+1$ according to whether $\{-\bar x\}\times\R$ or $\{\bar x\}\times\R$ is hit (and regardless of whether the threshold is hit from the left or from the right). The trajectory further travels along system $k$ until it reaches one of the switching lines again, when the same rule applies, see Figure~\ref{fig0}.
In this paper we only work with forward solutions that intersect the switching lines transversally.
 
 \vskip0.2cm

\noindent The paper deals with systems (\ref{fg})-(\ref{RL}) whose limit cycle shrinks to the origin when $\bar x\to 0
^+$. The latter can only happen when the vector fields $(f^-,g^-)$ and $(f^+,g^+)$ are opposite one another at $0$, i.e. when
\begin{equation}\label{sw_eq}
   \lambda \left(\begin{array}{c} f^-(0)  \\ g^-(0)\end{array}\right)+(1-\lambda)\left(\begin{array}{c}f^+(0) \\ g^+(0)\end{array}\right)=0,\quad\mbox{for some}\ \lambda\in[0,1],
\end{equation}
or, in other words, when $0$ is a switched equilibrium of the Filippov system
\begin{equation}\label{np}
  \left(\begin{array}{c}
      \dot x(t)\\ \dot y(t)\end{array}\right)=\left\{\begin{array}{l}\left(\begin{array}{l}f^-(x(t),y(t))\\ g^-(x(t),y(t))\end{array}\right), \quad{\rm if}\ x(t)<0, \\
\left(\begin{array}{l}f^+(x(t),y(t))\\ g^+(x(t),y(t))\end{array}\right), \quad{\rm if}\ x(t)>0.\end{array}\right.
\end{equation}
The fold-fold case where $f^-(0)=f^+(0)=0$ was considered in Makarenkov \cite{siam} in the context of an application to anti-lock braking systems. The present paper addresses the transversal case
$$
   f^-(0)f^+(0)\not=0
$$
as it is appears in the models of power converters. If the transversality condition holds then the switching threshold 
\begin{equation}\label{L}
\{(x,y)\in\mathbb{R}^2:x=0\}
\end{equation}
 is either an escaping or a sliding region of (\ref{np}) in the neighborhood of the origin. The case of escaping is not capable to produce limit cycles of (\ref{fg})-(\ref{RL}), so we assume the following stronger transversality condition 
\begin{equation}\label{positive}
   f^-(0)>0\quad{\rm and}\quad f^+(0)<0,
\end{equation}
which ensures sliding.

\vskip0.2cm

\noindent Limit cycles of switched system (\ref{fg})-(\ref{RL}) are of great importance in applied sciences. In addition to anti-lock braking systems in mechanical engineering (studied in \cite{siam,tan09,pas06}), these limit cycles are e.g. operating modes of intermittent therapy in medicine \cite{tanaka} and grazing management in ecology \cite{meza}. Any two-state control system (relay control system) in $\mathbb{R}^2$ is described by system (\ref{fg})-(\ref{RL}) as long as the sensor measures just one of the two state variables. The available results on the limit cycles are either for linear relay systems (e.g. Astrom \cite{astrom}, Andronov et al \cite[70, Ch. III, \S5]{andronov}) or for just existence (not stability) (e.g. Boiko \cite{boiko}). A particular application of the current paper will be in the field of  
switching power converters \cite{gupta,lu,schild} of electrical engineering. 

\vskip0.2cm

\noindent A dc-dc power converter that operates in a continuous current conduction mode switches between  states $k=-1$ and $k=1$ corresponding to the following two subsystems
\begin{eqnarray}
 &&\left(\begin{array}{c} (i_L)' \\  (u_C)' \end{array}\right)=\left(\begin{array}{cc} R_L/L & 0 \\ 0 & -1/(RC) \end{array}\right)\left(\begin{array}{c} i_L \\  u_C \end{array}\right)+\left(\begin{array}{c} V_s/L \\ 0 \end{array}\right),\qquad\qquad  {\rm if} \ \ k=-1,\label{dc1}\\
&& \left(\begin{array}{c} (i_L)' \\ (u_C)' \end{array}\right)=\left(\begin{array}{cc} R_L/L & -1/L \\ 1/C & -1/(RC) \end{array}\right)\left(\begin{array}{c} i_L \\  u_C \end{array}\right)+\left(\begin{array}{c} (V_s-V_d)/L \\ 0 \end{array}\right),\quad\ {\rm if}\ \ k=+1,\label{dc2}
\end{eqnarray}
where $i_L$ and $u_C$ are the inductor current and capacitor voltage respectively, $R_L,$ $L$, $R$, $C$, $V_s$ and $V_d$ are positive parameters of the circuit, see Schild et al \cite{schild}, 
Hu \cite{hu}, Gupta et al \cite {gupta}, Sreekumar et al \cite{sreekumar}
 for details. The switching rule of the form
  \begin{equation}\label{gensw}
     \begin{array}{ll}
        k:=+1, & {\rm if} \quad n(i_L,u_C)^T>c,\\
        k:=-1, & {\rm if}\quad n(i_L,u_C)^T<c,
     \end{array}
  \end{equation}
  where $n\in\mathbb{R}^2$ and $c\in\mathbb{R}$ are control parameters, makes system (\ref{dc1})-(\ref{gensw}) a Filippov system. Global stabilization of Filippov system (\ref{dc1})-(\ref{gensw}) near a switched equilibrium is addressed e.g. in Lu et al \cite{lu} using Lyapunov functions for convex combinations of linear systems, which strategy is originally suggested in Bolzern-Spinelli \cite{bolzern}. Our section~2 can be viewed as a local analog of this strategy. The Bolzern-Spinelli stabilization leads to sliding along the switching threshold $\{(i_L,u_C)\in\mathbb{R}^2:n(i_L,u_C)^T=c)\}$ and can also be achieved using the sliding mode control theory (see Fridman et al \cite{fridman}). When sliding is undesirable, different types of regularization of switching rule (\ref{gensw}) are used. We refer the reader to Lu et al \cite{lu}, Hu \cite{hu} and references therein for regularizations of (\ref{gensw}) which are capable to provide stabilization to a point. The main interest of the present paper is in the regularization of the form (termed {\it border-splitting} in Makarenkov-Lamb \cite{survey})
 \begin{equation}\label{genswreg}
     \begin{array}{ll}
        k:=+1, & {\rm if} \quad n(i_L,u_C)^T>c+\eps,\\
        k:=-1, & {\rm if}\quad n(i_L,u_C)^T<c-\eps,
     \end{array}
  \end{equation}
  which leads to a limit cycle. This type of regularization is used for dc-dc power converters e.g. in Schild et al \cite{schild}, where the existence of a cycle is discussed using the phase plane analysis. The present paper offers a bifurcation approach to the existence and orbital stability of a limit cycle in the system of (\ref{dc1})-(\ref{dc2}) and (\ref{genswreg}), which, in particular, gives new asymptotic formulas to assess properties of the limit cycle.

\vskip0.2cm

\noindent Perhaps surprisingly, there are several significant results on bifurcations of limit cycles from a switched equalibrium of (\ref{np}) under smooth perturbations of the vector fields of (\ref{np}) (see Guardia et al \cite{teixeira}, Kuznetsov et al \cite{kuznecov}, Simpson-Meiss \cite{simpson}, Zou et al \cite{zou}), but almost no research on bifurcation caused by splitting the switching boundary (also known as hysteresis perturbation). The respective results on bifurcations of limit cycles from a fold-fold singularity were recently obtained by Simpson \cite{simpson2014} and Makarenkov \cite{siam}, but border-splitting bifurcations from a switched equilibrium are missed in the literature. The current paper makes a significant progress in filling in this gap. 

\vskip0.2cm

\noindent Splitting of the switching manifold (also known as hysteresis switching or hybrid feedback switching in control, see Liberzon \cite{liberzon}) is an alternative to  piecewise linear continuous approximation (Sotomayor-Teixeira \cite{soto}) and to the singular perturbation approach (see Llibre-da Silva-Teixeira \cite{llibre}). Alike the regularization of \cite{llibre,soto}, border-splitting is supposed to unveil the actual behavior that an equilibrium of an idealized Filippov system will exhibit in the respective physical process. In other words, the dynamics coming from a border-splitting bifurcation is not always created on purpose (as in control applications discussed above), but can also be a universal consequence of switches. See e.g. Wojewoda et al \cite{woje} where the occurrence of hysteresis switching in dry friction oscillators is evidenced in experiments, or Cao-Chen \cite{actuators} where the effect of hysteresis is discussed in the context of electric actuators. 

\vskip0.2cm

\noindent The paper is organized as follows. Section 2 formulates the classical condition about the stability of the origin for Filippov system (\ref{np}). Section~3 shows (Theorem~\ref{thm1}) that if this classical condition holds for (\ref{np}), then an attracting limit cycle of the respective switched system (\ref{fg})-(\ref{RL}) emerges from the origin as the regularization parameter $\bar x$ crosses zero.  

\vskip0.2cm

\noindent The proof of Theorem~\ref{thm1} can serve as a guideline to a rigorous proof of bifurcation of limit cycles in piecewise-smooth systems, which is a good addition to the available textbooks in the field (see e.g. di Bernardo et al \cite{textbook}).
Section~4 applies the result to the dc-dc power converter from  Schild et al \cite{schild} and design a switching rule that leads to a limit cycle near a given reference value of the voltage variable.

\section{Local stability of the equilibrium of the Filippov system}

\noindent According to Filippov (see \cite[p.~218]{filippov}), the differential equation of sliding for (\ref{np}) along the threshold (\ref{L}) is given by
\begin{equation}\label{sliding}
   \dot y=\dfrac{f^+(0,y)g^-(0,y)-f^-(0,y)g^+(0,y)}{f^+(0,y)-f^-(0,y)}
\end{equation}
and the condition (\ref{sw_eq}) for the origin to be a switched equilibrium can be rewritten as
\begin{equation}\label{eq_cond}
   f^+(0)g^-(0)-f^-(0)g^+(0)=0.
\end{equation}
Computing the derivative of the right-hand-side of (\ref{sliding}) at the origin we conclude that the origin is an asymptotically stable equilibrium of (\ref{sliding}), if
\begin{equation}\label{stab}
   f^+{}'_y(0)g^-(0)+f^+(0)g^-{}'_y(0)-f^-{}'_y(0)g^+(0)-f^-(0)g^+{}'_y(0)>0.
\end{equation}

\section{The main result} 

\noindent In this section we prove that local stability of the origin for Filippov system (\ref{np}) always implies the existence of a stable limit cycle near the origin for the respective switched system (\ref{fg})-(\ref{RL}) when the regularization parameter $\bar x$ is sufficiently small.

\begin{thm}\label{thm1} Let $f$ and $g$ be $C^4$-functions. Assume that, for $\bar x=0$, the origin is an asymptotically stable equilibrium of Filippov system (\ref{np}), i.e.  conditions (\ref{positive}), (\ref{eq_cond}), and  (\ref{stab}) hold. Then, for all $\bar x>0$ sufficiently small, switched system (\ref{fg})-(\ref{RL}) possesses an orbitally stable limit cycle   which shrinks to the origin as $\bar x\to 0.$ The period $T(\bar x)$  of the cycle admits an estimate
\begin{equation}\label{Tbarx}
T(\bar x)=\left(-\frac{2}{f^+(0)}+\frac{2}{f^-(0)}\right)\bar x+O\left(\bar x^2\right).
\end{equation}
\end{thm}


\noindent {\bf Proof.} {\bf Step 1:} {\it The time map.} Let $\left(\begin{array}{c} X^k \\ Y^k \end{array}\right)(t,x,y)$ be the solution $t\mapsto(x(t),y(t))^T$ of system (\ref{fg}) with the initial condition $(x(0),y(0))=(x,y).$  Let $T^k(x,y)$ be the solution of the equation
$$
   X^k(t,x,y)=-x.
$$
Referring to Fig.~(\ref{fig0}b), $T^+(A)$ is the time that the solution of system (\ref{fg}) needs to go from point $A$ to point $B$ of $\mathbb{R}^2.$ Analogously, $T^-(B)$ is the time that the solution needs to go from point $B$ to point $C$. 
Such a function exists by the Implicit Function Theorem applied to
$$
  F(t,x,y)=X^k(t,x,y)+x,
$$
since $F'_t(0)=X^k{}'_t(0)=f^k(0)\not=0$ by (\ref{positive}). Expanding $T^k(x,y)$ in Taylor series we get
\begin{equation}\label{Ti}
   T^i(x,y)=T^i{}'_x(0)x+T^i{}'_x{}'_x(0)x^2+T^i{}'_x{}'_y(0)xy+O((x,y)^3),
\end{equation}
where we use that 
\begin{equation}\label{00}
T^i{}'_y(0)=T^i{}'_y{}'_y(0)=0,
\end{equation} which comes by taking the derivatives of 
\begin{equation}\label{Xi}
   X^i(T^i(x,y),x,y)=-x
\end{equation}
with respect to $y$ and using (\ref{positive}).

\vskip0.2cm
\noindent {\bf Step 2:} {\it Point transformation maps $P^i_x$ from $\{x\}\times\mathbb{R}$ to $\{-x\}\times\mathbb{R}$ and the Poincare map $P_x$ from $\{x\}\times\mathbb{R}$ to $\{x\}\times\mathbb{R}$.} Since
$$
   \frac{\partial}{\partial y}Y^k(T^k(x,y),x,y)=Y^k{}'_t(0)T^k{}'_y(0)+Y^k{}'_y(0)=1,
$$
we get the following expansion for the point transformation  $P_x^k(y)$ acting from $\{x\}\times\mathbb{R}$ to $\{-x\}\times\mathbb{R}$
 \begin{equation}\label{Px}
    P_x^k(y)=Y^k(T^k(x,y),x,y)=A^kx+y+a^kx^2+b^kxy+O((x,y)^3),\quad k=-1,+1,
 \end{equation}
 where we use that
 $$
    \left.\frac{\partial^2}{\partial y^2}Y^k(T^k(x,y),x,y)\right|_{(x,y)=0}=Y^k{}'_t{}'_t(0)\left(T^k{}'_y(0)\right)^2+2Y^k{}'_t{}'_y(0)T^k{}'_y(0)+Y^k{}'_y{}'_y(0)+Y^k{}'_t(0)T^k{}'_y{}'_y(0)=0
 $$
 by (\ref{00}). The  constants $a^i,b^i$ will be computed when necessary. Right now we only need a formula for $A^k$, so we compute 
 $$
   A^k=\left.\frac{\partial}{\partial x}Y^k\left(T^k(x,y),x,y\right)\right|_{(x,y)=0}=Y^k{}'_t(0)T^k{}'_x(0)+Y^k{}'_x(0)=g^k(0)T^k{}'_x(0).
 $$
To find $T^k{}'_x(0)$ we compute the derivative of (\ref{Xi}) with respect to $x$ and obtain 
\begin{equation}\label{Tx}
    X^k{}'_t(0)T^k{}'_x(0)+X^k{}'_x(0)=-1,\quad{\rm i.e.}\quad f^k(0)T^k{}'_x(0)+1=-1,\quad{\rm or}\quad T^k{}'_x(0)=-\frac{2}{f^k(0)},
\end{equation}
and 
$$
  A^k=-2\frac{g^k(0)}{f^k(0)}.
$$
Based on (\ref{eq_cond}) we now conclude that 
$$
   A^+-A^-=0
$$
and so,
the Poincare map $P_x$ takes the form 
\begin{equation}\label{Pmap}
   P_x(y)=P_{-x}^-(P^+_{x}(y))=y+ax^2+bxy+O((x,y)^3),
\end{equation}
whose constants $a$ and $b$ will be computed later, if required.
 
 \vskip0.2cm
 
\noindent {\bf Step 3:} {\it The existence of a fixed point.} Solving
$$
   P_x(y)=y
$$
one gets
\begin{equation}\label{of}
   ax^2+bxy+O\left((x,y)^3\right)=0\quad {\rm or}\quad ax+by+\frac{1}{x}O\left((x,y)^3\right)=0.
\end{equation}
Since $T^k(x,y)=T^k(x,y)-T^k(0,y)=\left(\int_0^1T^k{}'_x(\tau x,y)d\tau\right) x$ and $Y^k{}'_y{}'_y{}'_y(0,x,y)=0$, the $y^3$-term of $O\left((x,y)^3\right)$ of (\ref{of}) is of the form $h(x,y)x$, where $h$ is a $C^1$-function. As a consequence, in formula (\ref{of}) one has
$$
\frac{1}{x}O\left((x,y)^3\right)=O\left((x,y)^2\right).
$$
Therefore,  for all $x>0$ sufficiently small the unique fixed point $y(x)$ of $P_x$ that converges to the origin as $x\to 0$ is given by
\begin{equation}\label{relations}
   y(x)=-\frac{a}{b}x+O(x^2).
\end{equation}

\vskip0.2cm

\noindent {\bf Step 4:} {\it Verifying that $P_x^+$ maps $ y(x)$ from the past to the future. Verifying that $P_{-x}^-$ maps $P_x^+(y(x))$ from the past to the future.} By  (\ref{relations}) we have
\begin{equation}\label{relations1}
   P^+_x(y(x))=O(x).
\end{equation}
Plugging (\ref{relations}) and (\ref{relations1}) to (\ref{Ti}) and using (\ref{Tx}), we conclude that 
$$
   T^+(x,y(x))=-\frac{2}{f^+(0)}x+O(x^2)\quad {\rm and}\quad T^-(-x,P^+_x(y(x)))=\frac{2}{f^-(0)}x+O(x^2),
$$
i.e. both the time maps are positive by (\ref{positive}), for $x>0$ sufficiently small.

\vskip0.2cm

\noindent {\bf Step 5:} {\it Stability of the fixed point.} Computing the derivative of (\ref{Pmap}) yields
$$
  (P_x)'(y(x))=1+bx+\left.\frac{\partial}{\partial y}O\left((x,y)^3\right)\right|_{y=y(x)},
$$
i.e. $y(x)$ is a stable fixed point of $P_x$, if
$$
  b<0.
$$
\noindent {\bf Step 6:} {\it Computing $b$.} By definition,
$$
   b=b^+-b^-,
$$
where, using (\ref{00}), (\ref{Tx}), and $Y^k{}'_x{}'_y(0)=0,$
\begin{eqnarray*}
   b^k&=&\left.\frac{\partial^2}{\partial x\partial y}Y^k\left(T^k(x,y),x,y\right)\right|_{(x,y)=(0,0)}=\\
   &=&Y^k{}'_t{}'_t(0)T^k{}'_x(0)T^k{}'_y(0)+Y^k{}'_t{}'_y(0) T^k{}'_x(0)+Y^k{}'_t(0)T^k{}'_x{}'_y(0)+Y^k{}'_t{}'_x(0) T^k{}'_y(0)+Y^k{}'_x{}'_y(0)=\\
   &=&-Y^k{}'_t{}'_y(0) \frac{2}{f^k(0)}+Y^k{}'_t(0)T^k{}'_x{}'_y(0)\ =\ -g^k{}'_y(0)\frac{2}{f^k(0)}+g^k(0)T^k{}'_x{}'_y(0).
\end{eqnarray*}
Taking the derivative of (\ref{Xi}) with respect to $x$ and $y$ we get
$$
   X^k{}'_t{}'_t(0)T^k{}'_x(0)T^k{}'_y(0)+X^k{}'_t{}'_y(0)T^k{}'_x(0)+X^k{}'_t(0)T^k{}'_x{}'_y(0)+X^k{}'_t{}'_x(0)T^k{}'_y(0)+X^k{}'_x{}'_y(0)=0
$$
or
$$
  -f^k{}'_y(0)\frac{2}{f^k(0)}+f^k(0)T^k{}'_x{}'_y(0)=0,
$$
which allows to get
$$
   b^k= -g^k{}'_y(0)\frac{2}{f^k(0)}+2g^k(0)\frac{f^k{}'_y(0)}{f^k(0)^2}.
$$
Using (\ref{eq_cond}) we finally get
$$
     b=\frac{1}{f^+(0)f^-(0)}\left( f^+{}'_y(0)g^-(0)+f^+(0)g^-{}'_y(0)-f^-{}'_y(0)g^+(0)-f^-(0)g^+{}'_y(0)\right),
$$
which is negative by (\ref{positive}) and (\ref{stab}). 

\vskip0.2cm

\noindent The proof of the theorem is complete.\qed

\section{Application to a dc-dc power converter}

\noindent Considering the parameters 
\begin{equation}\label{XX}
   R_L=0.25\ \Omega,\ L=1\ mH,\ R=50\ \Omega,\ C=20.5\ \mu F,\ V_s=12\ V,\ V_d=0.4\ V,
   \end{equation}
  Schild et al \cite{schild} offers the following hybrid control strategy
  \begin{equation}\label{cs}
     \begin{array}{ll}
        k:=+1, & {\rm if} \quad (0.946,0.324)(i_L,u_C)^T> 6.44,\\
        k:=-1, & {\rm if}\quad (0.876,0.482)(i_L,u_C)^T< 9.2
     \end{array}
  \end{equation}
   to drive switched system of (\ref{dc1})-(\ref{dc2}) and (\ref{cs}) to a limit cycle around $V_{ref}=18\ V.$ We now show that a similar achievement can be obtained by applying Theorem~\ref{thm1}. Indeed, let us consider the switching rule
  \begin{equation}\label{cs1}
     \begin{array}{ll}
        k:=+1, & {\rm if} \quad n(i_L,u_C)^T> c+\eps,\\
        k:=-1, & {\rm if}\quad n(i_L,u_C)^T< c-\eps,
     \end{array}\quad{\rm where}\quad n=\left(\begin{array}{c} 0.91\\ 0.415\end{array}\right),
  \end{equation}
  where the normal vector $n$ is chosen to be roughly parallel to the two normal vectors of the switching rule (\ref{cs}). We will now first compute $c$ that gives $u_{ref}=18 V$ for the switched equilibrium of the respective Filippov system (\ref{np}) and then apply Theorem~\ref{thm1} in order to show the existence of an attracting limit cycle $t\mapsto(i_{L,\eps},u_{C,\eps})(t)$, 
  \begin{equation}\label{conv1}
  u_{C,\eps}(t)\to u_{ref},\quad {\rm as}\quad \eps\to 0^+,
  \end{equation} 
  to the switched system of (\ref{dc1})-(\ref{dc2}), (\ref{XX}), and (\ref{cs1}). 
   
\vskip0.2cm

\begin{figure}[h]\center
\includegraphics[scale=0.55]{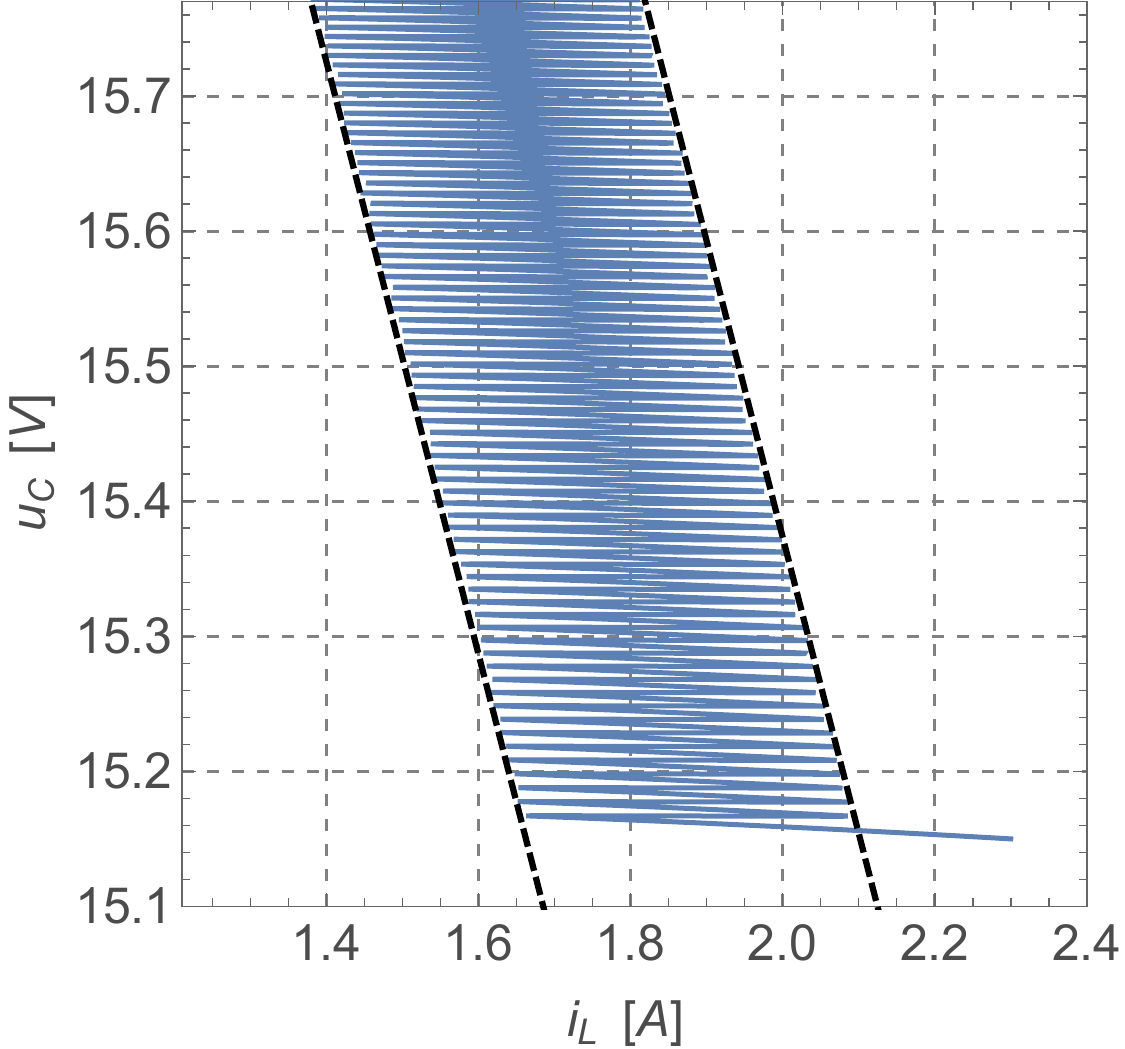}\ \ \   \includegraphics[scale=0.319]{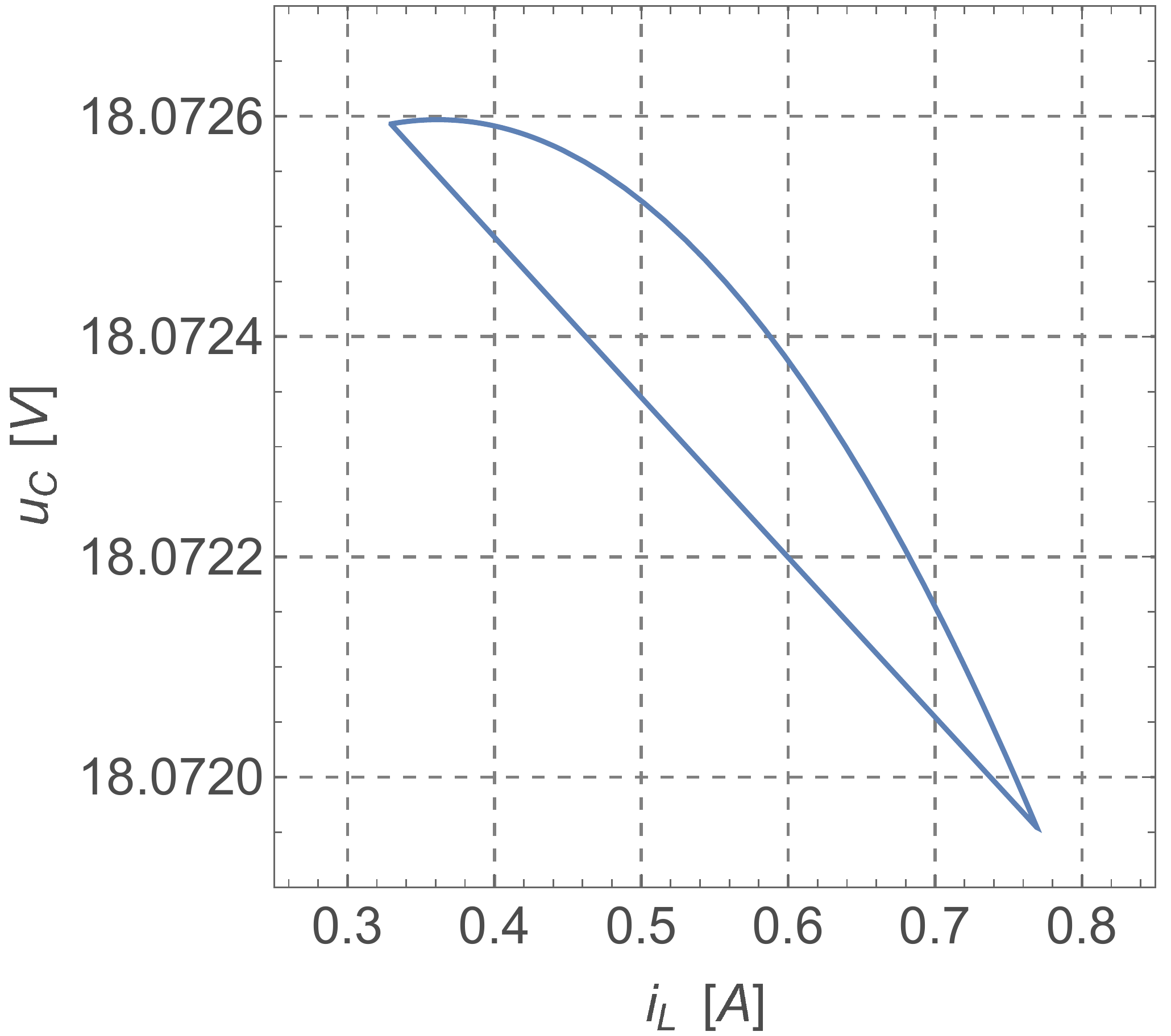}
\vskip-0.2cm
\caption{\footnotesize  Dynamics of the switched power converter model (\ref{dc1})-(\ref{dc2}) with  the parameters (\ref{XX}) and the switched control (\ref{cs1}) with $c=8$ and $\eps=0.2.$ The left and right dashed lines in each of the two figures are the two switching thresholds of (\ref{cs1}).  Left: The solution with the initial condition $(i_L(0),u_C(0))=(2.3,15.15)$. Right: The attractor of this solution (the attracting limit cycle of period $T\approx 0.1055$).} \label{fig1}
\end{figure}

\noindent The change of the variables
\begin{equation}\label{changev}
  \left(\begin{array}{c} i_L\\ u_C\end{array}\right)=\left(\begin{array}{cc} n_1 & -n_2 \\ n_2 & n_1\end{array}\right)\left(\begin{array}{c} x\\ y\end{array}\right)
\end{equation}   
transforms switched system (\ref{dc1})-(\ref{dc2}), and (\ref{cs1}) to the form
\begin{eqnarray}
\hskip-1.1cm  \left(\begin{array}{c}
    \dot x\\ \dot y\end{array}\right)&=&\left(\begin{array}{cc} n_1 & n_2 \\
    -n_2 & n_1\end{array}\right)\left[\left(\begin{array}{cc}
    -\frac{R_L}{L} & -\frac{k+1}{2}\cdot\frac{1}{L}\\
    \frac{k+1}{2}\cdot \frac{1}{C}&-\frac{1}{RC}
    \end{array}\right)\left(\begin{array}{cc} n_1 & -n_2 \\ n_2 & n_1\end{array}\right)\left(\begin{array}{c} x \\ y \end{array}\right)+ \left(\begin{array}{c} \frac{V_s}{L}-\frac{k+1}{2}\cdot \frac{V_d}{L} \\ 0\end{array}\right)\right], \label{ss1}\\
    k&:=&\left\{\begin{array}{ll}+1,& \quad{\rm if}\quad x>(c+\eps)/\|n\|^2,\\
    -1,& \quad {\rm if}\quad x<(c-\eps)/\|n\|^2.
    \end{array}\right. \label{ss2}
\end{eqnarray}  
With the parameters (\ref{XX}), switched system (\ref{ss1})-(\ref{ss2}) reads as 
\begin{eqnarray}
  \left(\begin{array}{c}
  \dot x\\
  \dot y
  \end{array}\right)&=&\left(\begin{array}{cc}
      -0.566 & -0.742\\
      0.307 & 0.315
      \end{array}\right)\left(\begin{array}{c}
   x\\
   y
  \end{array}\right)+\left(\begin{array}{c}
  10.556\\
  -4.814
  \end{array}\right),\quad {\rm if}\quad k=+1,\label{num1} \\
  \left(\begin{array}{c}
  \dot x\\
  \dot y
  \end{array}\right)&=&\left(\begin{array}{cc}
      -0.207 & 0.094\\
      0.094 & -0.044
      \end{array}\right)\left(\begin{array}{c}
   x\\
   y
  \end{array}\right)+\left(\begin{array}{c}
  10.92\\
  -4.98
  \end{array}\right),\quad {\rm if}\quad k=-1,\label{num2} \\
     k&:=&\left\{\begin{array}{ll}+1,& \quad{\rm if}\quad x>c+\eps,\\
    -1,& \quad {\rm if}\quad x<c-\eps.
    \end{array}\right. \label{num3}
\end{eqnarray}
The equilibrium $y$ of the respective Filippov system (\ref{sliding}) (that corresponds to $\eps=0$) is given by
\begin{equation}\label{yeq}
   -0.533c+0.243y+0.01c^2-0.008cy+0.003y^2=0.
\end{equation}
Based on (\ref{conv1})-(\ref{changev}), we further conclude that  the equilibrium $y$ of (\ref{yeq}) needs to satisfy  
\begin{equation}\label{yeq1}
  u_{ref}=n_2c+n_1y,\quad {\rm i.e.}\quad 18=0.415c+0.91y.
\end{equation}
Solving system (\ref{yeq})-(\ref{yeq1}) one gets 
$$
   c=7.976\approx 8,
$$
which needs to be used in the hybrid control rule (\ref{cs1}) to ensure property 
   (\ref{conv1}) with $V_{ref}=18.$ Numeric simulations of Fig.~\ref{fig1} support this conclusion and look similar to the respective simulations in \cite[Fig.~4]{schild} that use a slightly more complex switching rule compared to (\ref{gensw}).

\vskip0.2cm

\noindent Applying formula (\ref{Tbarx}) we get the following approximation for the period of the limit cycle of system (\ref{dc1})-(\ref{dc2}) with parameters (\ref{XX}) and switched rule (\ref{cs1}) with $c=8$  and $\eps=0.2$:
$$
   T(0.2)=0.1043,
$$
which agrees with the period of $0.1055$ computed numerically (see caption of Fig.~\ref{fig1}) quite well.

\section{Conclusion}  \noindent We proved (Theorem~\ref{thm1}) that a regular asymptotically stable equilibrium of the Filippov equation of sliding motion (known as {\it switched equilibrium} in control) yields an orbitally stable limit cycle upon a hysteresis type switched regularization  (known as {\it hybrid feedback switching} in control).  The research is motivated by applications in power electronics where switched equilibrium is a standard operating mode of switched dc-dc power converters. Regularization is usually introduced in the control strategy of the dc-dc converters in order to avoid sliding modes that occur for the most directed switched control strategies (based on Lyapunov functions of stable convex combinations). 

\vskip0.2cm

\noindent The proposed formula (\ref{Tbarx}) of Theorem~\ref{thm1} for the period of the limit cycle says that the period of the limit cycle grows linear with respect to the regularization parameter. This is similar to the linear growth  observed the other nonsmooth analogues of Hopf bifurcations known in the literature (as opposed to the square-root growth typical for the smooth Hopf bifurcation), see Simpson-Meiss \cite{aspects}. Our  numeric simulation with a particular dc-dc power converter model (Section~5) discovered a very good accuracy of formula (\ref{Tbarx}). This formula can therefore be used as a tool to control the period of the limit cycle in dc-dc power converters with hysteresis switching rules.

\vskip0.2cm

\noindent When proving Theorem~\ref{thm1} we established  the fact that the initial condition $(\bar x,y(\bar x))$ of the limit cycle of switched system (\ref{fg})-(\ref{RL}) is given by the formula
$$
   y(\bar x)=-\frac{a}{b}\bar x+O\left(\bar x^2\right),
$$
see formula (\ref{relations}), where we never computed the value $a$ as we didn't need it for the required conclusions. However, one can go further and compute a closed-form formula for $a$ and get an asymptotic estimate for the initial condition $(\bar x,y(\bar x))$.  Such an asymptotic estimate can be further used in order to estimate the angle of the corners of the limit cycle (think of the corners of the limit cycle of Fig.~\ref{fig1}right). The later will provide a tool  to control the thickness of the limit cycle.

\section*{Compliance with Ethical Standards}

\noindent {\bf Acknowledgment:} The study was funded by the National Science Foundation (grant 
CMMI-1436856). 

\vskip0.2cm

\noindent {\bf Conflict of Interest:} The author has no conflict of interest.

\bibliographystyle{plain}

\end{document}